\documentclass{amsart}
\usepackage{amssymb,latexsym,graphics,amsfonts}

\swapnumbers
 \theoremstyle{plain}

 \theoremstyle{definition}
 
 \theoremstyle{remark}
 
 \newcommand{\cal}[1]{\mathcal{#1}}

\begin{document}
\title{Homomorphisms, amenability and weak amenability of Banach
algebras}
\author{M. Eshaghi Gordji}
\address{Department of Mathematics,
University of Semnan, Semnan, Iran} \email{maj\_ess@Yahoo.com and
meshaghi@semnan.ac.ir}
\keywords{Derivation , Homomorphism} \subjclass[2000]{46HXX}
\dedicatory{}
\smallskip
\begin{abstract}
In this paper we find some necessary and sufficient conditions for
a Banach algebra to be amenable or weakly amenable, by applying
the homomorphisms on Banach algebras.
\end{abstract}
\maketitle
\section{Introduction}
Let $\cal A$ be a Banach algebra and let $X$ be a Banach $\cal
A$-bimodule.
 Then $X^*$ is a Banach $\cal A$-bimodule if for each $a\in \cal A$
, $x\in X$ and $x^*\in X^*$ we define \[\langle
x,ax^*\rangle=\langle xa,x^*\rangle, \hspace{1cm} \langle
x,x^*a\rangle=\langle ax,x^*\rangle.\] Let $\varphi:\cal
A\longrightarrow\cal B$ be a Banach algebra homomorphism, then
$\cal B$ is a $\cal A$-bimodule by the following module actions
$$a.b=\varphi(a)b , \hspace{0.7cm} b.a=b\varphi(a) \hspace{1.5cm} (a\in
\cal A, b\in \cal B).$$ We denote $\cal B_{\varphi}$ the above $\cal
A$-bimodule. For a Banach algebra ${\cal A}$, ${\cal A}^{**}$ with
the first Arens product is a Banach algebra. Let $X$ be a Banach
$\cal A-$module, we can extend the actions of ${\cal A}$ on $X$ to
actions of ${\cal A}^{**}$  on $X^{**}$ via
\[a''.x'' =w^{*}\textup{-}\lim_{i} \lim_{j} a_{i} \, x_{j}
\]
and
\[ x''.a'' =w^{*}\textup{-}\lim_{j} \lim_{i} x_{j} \, a_{i},\]
where $a'' =w^{*}\textup{-}\lim_{i} a_{i}$, \  $x''
=w^{*}\textup{-}\lim_{j} x_{j}$.

If $X$ is a Banach $\cal A$-bimodule then a derivation from $\cal A$
into $X$ is a continuous linear map $D$, such that for every $a,b
\in \cal A,$ $D(ab)=D(a).b+a.D(b).$ If $x\in X,$ and we define
$\delta_x:\cal A\longrightarrow X$ by $\delta_x(a)=a.x-x.a\hspace
{0.5cm}(a\in \cal A),$ then $\delta_x$ is a derivation, derivations
of this form are called inner derivations. A Banach algebra $\cal A$
is amenable if every derivation from $\cal A$ into each dual $\cal
A$-bimodule is inner; i.e. $H^1(\cal A,X^*)=\{o\}$, this definition
was introduced by B. E. Johnson in [7]. $\cal A$ is weakly amenable
if $H^1(\cal A,\cal A^*)=\{0\},$ where $H^1(\cal A,\cal A^*)$ is the
first cohomology group from $\cal A$ with coefficients in $\cal
A^*$. Bade, Curtis and Dales have introduced the concept of weak
amenability for commutative Banach algebras [1]. In this paper we
show that for amenability of Banach algebra $\cal A$, it is enough
to show that for every Banach algebra $\cal B$ and every injective
homomorphism $\varphi:\cal A\longrightarrow\cal B$, $H^1(\cal
A,{\cal B_{\varphi}}^*)=\{o\}.$ So we introduce two new notations in
amenability of Banach algebras and we related this notations to weak
amenability.

\section{Amenability}
Let $\cal A$  be a Banach algebra and $X$ be a Banach $\cal
A$-bimodule, then $X\oplus_1 {\cal A}$ is a Banach space, with the
following norm
$$\|(x,a)\|=\|x\|+\|a\|\hspace{1.5cm} (a\in{\cal A}~,~x\in X).$$
So $X\oplus_1 {\cal A}$ is a Banach algebra with the product
$$(x_1,a_1)(x_2,a_2)=(x_1\cdot a_2+a_1\cdot x_2,a_1a_2)~.$$
$X\oplus_1 {\cal A}$ is called a module extension Banach algebra.
It is easy to show that $({X\oplus_1 {\cal A}})^*=X^*\oplus {\cal
A}^*$, where this sum is $\cal A$-bimodule $l_\infty-$sum. In this
section we use module extension Banach algebras to finding an easy
equivalent condition for amenability of a Banach algebra.
\paragraph{\bf Theorem 2.1.}
Let $\cal A$ be a Banach algebra. Then the following assertions
are
equivalent:\\
(i) $\cal A$ is amenable.\\
(ii) For every Banach algebra $\cal B$ and every homomorphism
$\varphi:\cal A\longrightarrow\cal B$, $H^1(\cal A,\cal
{B_{\varphi}}^*)=\{o\}.$\\
(iii) For every Banach algebra $\cal B$ and every injective
homomorphism $\varphi:\cal A\longrightarrow\cal B$, $H^1(\cal
A,\cal
{B_{\varphi}}^*)=\{o\}.$\\
(iv) For every Banach algebra $\cal B$ and every injective
homomorphism $\varphi:\cal A\longrightarrow\cal B$, if
$d_{\varphi}:\cal A\longrightarrow {\cal B_{\varphi}}^{*}$ is a
(bounded) derivation satisfies
$$\langle d_{\varphi}(a),\varphi
(b)\rangle+\langle d_{\varphi}(b) ,\varphi (a)\rangle=0 \hspace{1cm}
(a,b\in \cal A),$$ then $d_{\varphi}$ is inner derivation.\\
(v) For every Banach algebra $\cal B$ and every injective
homomorphism $\varphi:\cal A\longrightarrow\cal B$, $H^1(\cal A,\cal
{B_{\varphi}}^{**})=\{o\}.$
\paragraph{\bf Proof.} The proofs of $(i)\Longrightarrow (ii)$, $(i)\Longrightarrow (v)$,  $(ii)\Longrightarrow
(iii)$ and $(iii)\Longrightarrow (iv)$ are immediate. We show that
$(iv)\Longrightarrow (i)$ and $(v)\Longrightarrow (i)$ hold. Suppose
that (iv) holds and let $X$ be a Banach $\cal A$-bimodule and
$D:\cal A\longrightarrow X^*$ be a derivation. As above, we khnow
that  $X\oplus_1 {\cal A}$ is a Banach algebra and obviously the map
$$\varphi:  a\mapsto(o,a),\hspace{0.7cm} \cal A\longrightarrow
X\oplus_1 {\cal A}$$ is an injective Banach algebra homomorphism.
Then $H^1(\cal A, ({(X\oplus_1 {\cal A})_{\varphi}})^*)=\{o\}.$ We
define $D_1:\cal A\longrightarrow (X\oplus_1 {\cal A})^*$ by
$D_1(a)=(D(a),o)$. For $a,b \in \cal A$ we have
\begin{align*} D_1(ab)&=(D(ab),0)=(D(a)b+aD(b),0) \\
&=(D(a),0)(0,b)+(0,a)(D(b),0)\\
&=D_1(a)\varphi(b)+\varphi(a)D_1(b).
\end{align*}
Thus $D_1$ is a derivation from $\cal A$ into $({(X\oplus_1 {\cal
A})_{\varphi}})^*$. Also for every $a,b\in \cal A,$   we have
$$\langle D_1(a),\varphi (b)\rangle+\langle D_1(b) ,\varphi (a)\rangle=\langle (D(a),0),(0,b)\rangle+\langle (D(b),0),(0,a)\rangle=0.$$
Then $D_1$ is inner derivation. On the other word there exist
$a'\in \cal A^* , x'\in X^*$ such that $D_1=\delta_{(x',a')}.$ For
every $a\in \cal A$ we have
\begin{align*} (D(a),o)&=D_1(a)=\delta_{(x',a')}(a) \\
&=\varphi(a)(x',a')-(x',a')\varphi(a) \\
&= (0,a)(x',a')-(x',a')(0,a) \\
&=(ax'-x'a,aa'-a'a).
\end{align*}
Thus $D=\delta_{x'}$. So $\cal A$ is amenable. To prove
$(v)\Longrightarrow (i)$, let  $X$ be a Banach $\cal A$-bimodule and
let $D:\cal A\longrightarrow X^{**}$ be a derivation. If $\varphi:
\cal A\longrightarrow X\oplus_1 {\cal A}$ is the above injective
Banach algebra homomorphism, then it is easy to show that
$\varphi^{**}: \cal A^{**}\longrightarrow (X\oplus_1 {\cal A})^{**}$
the second transpose of $\varphi$ is a Banach algebra homomorphism
and that $({(X\oplus_1 {\cal A})_{\varphi}})^{**}\simeq
({X^{**}\oplus_1 {\cal A}^{**})_{\varphi^{**}}}$ as $\cal
A^{**}$-bimoduls. Then
$$H^1(\cal A, ({X^{**}\oplus_1 {\cal
A}^{**})_{\varphi^{**}}})=H^1(\cal A, ({(X\oplus_1 {\cal
A})_{\varphi}})^{**})=\{o\} \hspace{0.7cm} (1).$$ Now we define
$D_1:\cal A\longrightarrow X^{**}\oplus_1 {\cal A}^{**}$ by
$D_1(a)=(D(a),o)$. For $a,b \in \cal A$ we have
$$ D_1(ab)=D_1(a)\varphi^{**}(\hat{b})+\varphi^{**}(\hat{a})D_1(b). $$ Thus
$D_1$ is a derivation from $\cal A$ into $(X^{**}\oplus_1 {\cal
A}^{**})_{{\varphi}^{**}}$. By (1), $D_1$ is inner. Therefore there
exist $a''\in \cal A^{**} , x''\in X^{**}$ such that
$D_1=\delta_{(x'',a'')},$ and by a similar proof as above we can
show that $D$ is inner. Then we have  $H^1(\cal A,X^{**})=\{o\}$,
and by Proposition 2.8.59 of [2], $\cal A$ is
amenable.\hfill$\blacksquare~$

Let $\cal A$ has a bounded approximate identity, and let $X$ be an
essential Banach $\cal A$-bimodule, then it is easy to show that
$({X\oplus_1 {\cal A}})_\varphi$ is an essential Banach $\cal
A$-bimodule when $\varphi:\cal A\longrightarrow X\oplus_1 {\cal
A}$ defined by $\varphi (a)=(o,a).$ By a same technique as above
and by using Corollary 2.9.28 of [2], we have the following
\paragraph{\bf Theorem 2.2.}
Let $A$ be a Banach algebra with a bounded approximate identity.
Then $\cal A$ is amenable if and only if for every Banach algebra
$\cal B$ and every injective homomorphism $\varphi:\cal
A\longrightarrow\cal B$ in which $\cal B_\varphi$ is essential,
$H^1(\cal A,{\cal B_{\varphi}}^*)=\{o\}.$
\section{ Weak amenability}
In this section we fined the relationship between weak amenability
and homomorphisms of Banach algebras. First we introduce two new
notations of amenability of Banach algebras.
\paragraph{\bf Definition 3.1.} Let ${\cal A}$ be a Banach algebra.
Then\\
(i) $\cal A$ is supper weakly amenable if for every Banach algebra
$\cal B$ and every continuous homomorphism $\varphi:\cal A
\longrightarrow \cal B$, if $d_{\varphi}$ is a (bounded)
derivation from ${\cal A}$ into ${\cal B_{\varphi}}^{*}$, then the
following condition holds
$$\langle d_{\varphi}(a),\varphi (b)\rangle+\langle d_{\varphi}(b)
,\varphi (a)\rangle=0  \hspace{1cm} (a,b\in \cal A)
\hspace{1cm}(1).$$ (iii) $\cal A$ is semiweakly amenable if every
derivation $D:\cal A\longrightarrow \cal A^*,$ by the following
property
$$\langle D(a),b\rangle+\langle D(b),a\rangle=0  \hspace{1cm} (a,b\in
\cal A)  \hspace{1cm}(2),$$ is an inner derivation.
\paragraph{\bf Example 1.} Let $\Bbb T$ be the unit circle. We write $(\hat {f}(n):n\in \Bbb
Z)$ for the sequence of Fourier coefficients of a function $f\in
L^1(\Bbb T)$. For $\alpha \in (\frac{1}{2},1)$, let $\cal
A=lip_\alpha (\Bbb T)$, we define $D:\cal A\longrightarrow \cal
A^*$ by
$$\langle D(f), g\rangle=\sum n \hat{g}(n) \hat {f}(-n),
\hspace{1cm}(f,g \in \cal A).$$ $D$ is a non-inner derivation (see
[1]) and we have
$$\langle D(f), g\rangle +\langle D(g), f\rangle=0 \hspace{1cm} (f,g
\in \cal A). $$ Thus $\cal A$ is not semiweakly amenable.
\paragraph{\bf Theorem 3.2.} Let ${\cal A}$ be a supper weakly amenable
Banach
algebra. Then\\
(i) $\cal A$ is essential.\\
(ii) There are no no-zero continuous point derivations on $\cal
A.$
\paragraph{\bf Proof.}(i): Let $a_0\in \cal A- \bar{\cal A^2},$ then by Hahn-Banach theorem
there exists $f\in \cal A^*$ such that $\langle f,a_0 \rangle=1$
and $f(\bar{\cal A^2})=\{0\}.$ The mapping $D:a \longmapsto f(a)f,
\hspace{0.5cm} \cal A\longrightarrow \cal A^*$ is a derivation and
we have $\langle D(a_0),a_0\rangle+\langle D(a_0),a_0\rangle=2\neq
0$. Thus ${\cal A}$ is not supper weakly amenable. (ii): Let
$\varphi\in \Omega_{\cal A}$. If $\varphi=0$, then by (i), every
derivation from $\cal A$ into ${\Bbb C_{\varphi}}^*$ is zero. If
$\varphi\neq 0,$ and  $d_{\varphi}: \cal A\longrightarrow \Bbb
C_{\varphi}$ is a point derivation at $\varphi$, then by
Definition 3.1, for every $a \in \cal A$ we have $\langle
d_\varphi(a),\varphi (a) \rangle=d_\varphi(a)\varphi(a)=0.$
Therefor we have $d_\varphi | (\cal A\setminus M_\varphi)=0.$ Thus
$d_\varphi =0.$\hfill$\blacksquare~$
\paragraph{\bf Example 2.} Let $\cal A=\Bbb C$ by the product $ab=0, (a,b\in \Bbb C).$ Then
by Theorem 3.2 (i), $\cal A$ is not supper weakly amenable. But it
is easy to check that $\cal A$ is semiweakly amenable.
\paragraph{\bf Example 3.} Let $S$ be a discrete semigroup in which $S^2\neq S$,
then by Theorem 3.2 (i), $l^1(S)$ is not supper weakly amenable.
Let $S=\{t,0\}$ by products $t0=0t=t^2=0^2=0$, then $l^1(S)$ is
not supper weakly amenable but for every derivation
$D:l^1(S)\longrightarrow {l^1(S)}^*$ if $\langle
D(\delta_t),\delta_0\rangle +\langle
D(\delta_0),\delta_t\rangle=0$, then we have $D=0$. Thus $l^1(S)$
is semiweakly amenable.

 Now we fined an equivalent condition for weak amenability of Banach
algebras.
\paragraph{\bf Theorem 3.3.} Let ${\cal A}$ be a Banach algebra.
Then \\
(i) $\cal A$ is weakly amenable if and only if $\cal A$ is supper
weakly amenable and semiweakly amenable.\\
(ii) Let $\cal A$ be a unital Banach algebra then $\cal A$ is
supper weakly amenable if and only if for every derivation $D:\cal
A\longrightarrow \cal A^*$, and for every $a \in \cal A,$ we have
$\langle D(a),1\rangle=0.$
\paragraph{\bf Proof.} (i): Let $\cal A$ be weakly
amenable. Obviously $\cal A$ is semiweakly amenable. For Banach
algebra
 $\cal B$ and for (continuous) homomorphism $\varphi:\cal A
\longrightarrow \cal B$,  let $d_{\varphi}:{\cal
A}\longrightarrow{\cal B_{\varphi}}^{*}$ be a derivation. We
define
 $D=d_\varphi\otimes \varphi :{\cal A}\longrightarrow{\cal A}^{*}$ as
 follows
 $$\langle D(a),b \rangle=\langle d_{\varphi}(a),\varphi(b)\rangle
\hspace{1cm}(a,b \in \cal A) \hspace{1cm}(3).$$ For every $a,b,c
\in \cal A$, we have
\begin{align*} \langle D(ab),c\rangle &=\langle
d_{\varphi}(ab),\varphi(c)
\rangle\\
&=\langle d_{\varphi}(a)\varphi (b),\varphi(c) \rangle+\langle
\varphi (a)d_{\varphi}(b),\varphi(c) \rangle\\
&=\langle d_{\varphi}(a),\varphi (b)\varphi(c) \rangle+\langle
d_{\varphi}(b),\varphi(c)\varphi (a) \rangle\\
&=\langle d_{\varphi}(a),\varphi (bc) \rangle+\langle
d_{\varphi}(b),\varphi(ca) \rangle\\
&=\langle D(a),bc \rangle+\langle D(b),ca \rangle\\
&=\langle D(a)b+aD(b),c \rangle.
\end{align*}
Therefore $D$ is a derivation. Then there exists $f\in \cal A^*$
such that $D=\delta_f:\cal A\longrightarrow \cal A^*$. Thus for
every $a,b\in \cal A,$ we have
\begin{align*}\langle
D(a),b\rangle+\langle D(b),a\rangle&=\langle
\delta_f(a),b\rangle+\langle \delta_f(b),a\rangle\\
&= \langle af-fa,b\rangle+\langle bf-fb,a\rangle\\
&=0. \end{align*} So $\cal A$ is supper weakly amenable. The
converse is trivially since $id:\cal A\longrightarrow \cal A$ is a
homomorphism in which $\cal A^*={\cal A_{id}}^*.$ (ii): Let  for
every derivation $D:\cal A\longrightarrow \cal A^*$, and for every
$a \in \cal A,$ the equality $\langle D(a),1\rangle=0$ holds, and
let $\varphi:\cal A \longrightarrow \cal B$ be a Banach algebra
homomorphism. If $d_{\varphi}:{\cal A}\longrightarrow{\cal
B_{\varphi}}^{*}$ is a derivation, then
 $D=d_\varphi\otimes \varphi :{\cal A}\longrightarrow{\cal A}^{*}$
defined in the proof of (i),
 is a derivation and for every $a,b\in \cal A,$ we have
$$\langle d_{\varphi}(a),\varphi (b)\rangle+\langle
d_{\varphi}(b),\varphi(a) \rangle=\langle D(a),b\rangle+\langle
D(b),a \rangle=\langle D(ab),1\rangle=0.$$ The converse is
trivial. \hfill$\blacksquare~$
\paragraph{\bf Corollary 3.4.} (Theorem 2.8.63 of [2]) Let $\cal
A$ be a weakly amenable Banach algebra, then $\cal A$ is essential
and there are no non-zero, (continuous) point derivations on $\cal
A.$
\paragraph{\bf Corollary 3.5.} Let $G$ be a locally compact
topological group, Then $G$ is discrete if and only if $M(G)$ is
supper weakly amenable.
\paragraph{\bf Proof.} H. G. Dales,  F. Ghahramani and A. Ya.
Helmeskii [3] showed that $G$ is discrete if and only if there are
no nonzero point derivations on $M(G)$. By applying Theorems 3.2
(ii) and 3.3 (i), we conclude that $G$ is discrete if and only if
$M(G)$ is supper weakly amenable.\hfill$\blacksquare~$

 By the following Theorem we can show that the supper weak amenability is
different from the weak amenability and semiweak amenability.
\paragraph{\bf Theorem 3.6.} Let $\cal A$ be a
 supper weakly amenable Banach algebra, and let  $\theta:\cal A\longrightarrow\cal B$
be a continuous Banach algebra homomorphism with dense range. Then
$\cal B$ is supper weakly amenable.
\paragraph{\bf Proof.} Let $\varphi:\cal B
\longrightarrow \cal C$ be a Banach algebra homomorphism and let
$d_{\varphi}:{\cal B}\longrightarrow{\cal C_{\varphi}}^{*}$ be a
derivation. Then for every $a,b\in \cal A,$ we have
$$d_{\varphi}o\theta(ab)=d_{\varphi}o\theta(a)\varphi o\theta (b)+\varphi o \theta (a)
d_{\varphi}o\theta(b).$$ Therefore $d_{\varphi}o\theta$ is a
derivation from $\cal A$ into $({\cal C_{\varphi o \theta}})^{*}.$
Since $\cal A$ is  supper weakly amenable, then for every $a,b\in
\cal A,$ we have
$$\langle {d_{\varphi}o\theta}(a),\varphi o\theta(b)\rangle+\langle
{d_{\varphi}o\theta}(b),\varphi o\theta(a)\rangle=0.$$ Since
$\theta (\cal A)$ is dense in $\cal B,$ then for every $a',b'\in
\cal B,$
$$\langle {d_{\varphi}}(a'),\varphi (b')\rangle+\langle
{d_{\varphi}}(b'),\varphi (a')\rangle=0.$$ Thus $\cal B$ is supper
weakly amenable. \hfill$\blacksquare~$
\paragraph{\bf Corollary 3.7.} There exists a supper weakly amenable,
non-semiweakly amenable Banach algebra.
\paragraph{\bf Proof.} Let $E$ be Banach space without approximation
property and take $\cal A$ to be the nuclear algebra
 $E\hat{\otimes} E^{*}$ (see Definition 2.5.4 of [2]). The identification of
 $E\otimes E^{*}$ with $\cal F(E)$ extends to an  epimorphism
$R:E\hat{\otimes} E^{*}\longrightarrow \cal N(E)$ (see Theorem
2.5.3 of [2]). Set $K=ker R,$ then by Corollary 2.8.43 of [2],
$\cal A$ is biprojective and hence weakly amenable. If $dim K\geq
2$ then $K$ does not have trace extension property. So by
Proposition 2.8.65 (iii) of [2], $\cal N(E)=\frac{\cal A}{K}$ is
not weakly amenable. On the other hand by (i) of Theorem 3.3,
above, $\cal A$ is supper weakly amenable and by Theorem 3.6,
$\cal N(E)=\frac{\cal A}{K}$ is supper weakly amenable. Thus  by
Theorem 3.3, $\cal N(E)$ is supper weakly amenable, non-semiweakly
amenable Banach algebra.\hfill$\blacksquare~$

We finish this section with a Theorem about semiweak amenability of
unitization of Banach algebras, and its application to finding an
example of non-supper weakly amenable Banach algebra which its
unitization is supper weakly amenable.
\paragraph{\bf Theorem 3.8.} Let ${\cal A}$ be a Banach algebra. If
$\cal A^\sharp$ (the unitization of $\cal A$) is
 semiweakly amenable, then $\cal A$ is semiweakly amenable.
\paragraph{\bf Proof.} Let $D:\cal A
\longrightarrow \cal A^*$ be a derivation in which (2) holds. We
define $D^\sharp:\cal A^\sharp\longrightarrow {\cal A^\sharp}^*$ as
follows
$$\langle D^\sharp(a,c),(b,c') \rangle=\langle D(a),b \rangle \hspace{1cm}
(a,b\in \cal A, c,c'\in \Bbb C).$$ Then for every $a,b,d\in \cal A$
and  $c,c',c''\in \Bbb C,$ we have
\begin{align*}\langle
D^\sharp((a,c)(b,c')),(d,c'')\rangle&=\langle
D^\sharp(ab+cb+c'a,cc'),(d,c'')\rangle=\langle
D(ab+cb+c'a),d\rangle\\&=\langle
D(a)b+aD(b)+cD(b)+c'D(a),d\rangle\\
&= \langle D(a),bd+c''b+c'd\rangle
+\langle D(b),da+cd+c''a\rangle\\
&= \langle D^\sharp(a,c),(bd+c''b+c'd,c'c'')\rangle+\langle
D^\sharp(b,c'),(da+cd+c''a,cc'')\rangle\\
&=\langle (D^\sharp(a,c))(b,c'),(d,c'')\rangle+\langle
(a,c)(D^\sharp(b,c')),(d,c'')\rangle.
\end{align*}
Thus $D^\sharp$ is a derivation. So we have
$$\langle D^\sharp(a,c),(b,c') \rangle=\langle D(a),b \rangle=\langle D(b),a \rangle=\langle D^\sharp(b,c'),(a,c) \rangle \hspace{1cm}
(a,b\in \cal A, c,c'\in \Bbb C).$$
 $\cal A^\sharp$ is semiweakly amenable, then there is $u'\in
{\cal A^\sharp}^*$ such that $D^\sharp=\delta_{u'}.$ So we have
$D=\delta_{(u'\mid_{\cal A})}.$ Thus $\cal A$ is semiweakly
amenable.\hfill$\blacksquare~$

Let $\cal A$ be the augmentation ideal of $L^1(PS(2,\Bbb R)),$
then we know that $\cal A^\sharp$, is weakly amenable and that
$\cal A$ is not weakly amenable (see [8]). By above Theorem, $\cal
A$ is semiweakly amenable. So by Theorem 3.3, $\cal A^\sharp$ is
supper weakly amenable and $\cal A$ is not supper weakly amenable.
Thus we have the following.
\paragraph{\bf Corollary 3.9.} There exists a semiweakly amenable Banach algebra $\cal A$, Such that
$\cal A^\sharp$ is supper weakly amenable, and $\cal A$ is not
supper weakly amenable.
\section{ supper weak amenability of the second dual of Banach algebras}
Let $\cal A^{**}$ be the second dual of $\cal A$ with the first
Arens product. Then amenability of $\cal A^{**}$ implies the
amenability of $\cal A$ (see for example Proposition 2.8.59 of
[2]). So weak amenability of $\cal A^{**}$ implies the weak
amenability of $\cal A$ if  one of the following conditions holds
(see [4],
[5] and [6]).\\
(i) ${\cal A}$ is a left ideal in ${\cal A}^{**}$.\\
(ii) ${\cal A}$ is a dual Banach algebra.\\
(iii) ${\cal A}$ is Arens regular and every derivation from ${\cal
A}$ into its dual is weakly compact.\\
Similarly for supper weak amenability we have
\paragraph{\bf Theorem 4.1.} Let $\cal A$ be a Banach algebra with
one of the conditions (i), (ii) or (iii) as above. Let $\cal
A^{**}$ be supper weakly amenable, then $\cal A$ is supper weakly
amenable.
\paragraph{\bf Proof.} Let $D:\cal A\longrightarrow \cal A^*$ be a derivation,
then $D$ has an extension $\bar{D}:\cal A^{**}\longrightarrow
{({\cal A}^{**})}^*$ in which $\bar{D}$ is a derivation (see [4],
[5] and [6]). Since $\cal A^{**}$ is supper weakly amenable, then
for every $a,b \in \cal A,$ we have
$$\langle D(a),b \rangle+\langle D(b),a \rangle=\langle \bar {D}(\hat
{a}),\hat {b} \rangle+\langle \bar {D}(\hat {b}),\hat {a}
\rangle=0.$$ Now let $\varphi:\cal A \longrightarrow \cal B$ be a
Banach algebra homomorphism and let $d_{\varphi}:{\cal
A}\longrightarrow{\cal B_{\varphi}}^{*}$ be a derivation. Set
$D=d_\varphi \otimes \varphi$ (defined in the proof of Theorem
3.3). We have
$$\langle d_{\varphi}(a),\varphi (b)\rangle+\langle d_{\varphi}(b)
,\varphi (a)\rangle=\langle D(a),b \rangle+\langle D(b),a
\rangle=0.$$ Thus $\cal A$ is supper weakly
amenable.\hfill$\blacksquare~$
\paragraph{\bf Theorem 4.2.} Let $\cal A$ be a Banach algebra. Let
$\cal A^{**}$
be supper weakly amenable, then\\
(i) $\cal A$ is essential.\\
(ii) There are no no-zero continuous point derivations on $\cal
A.$
\paragraph{\bf Proof.} (i): By Theorem 3.2 (i), $\cal A^{**}$ is essential,
then we can show that $\cal A$ is essential (see for example
Proposition 2.1 of [5]). (ii): Let  $\varphi\in \Omega_{\cal A}$.
If $\varphi=0$, then by (i), every derivation from $\cal A$ into
${\Bbb C_{\varphi}}^*$ is zero. Now let $0\neq \varphi\in
\Omega_{\cal A}$, then it is easy to show that ${\varphi}''\in
\Omega_{{\cal A}^{**}}.$ We suppose that $d_{\varphi}: \cal
A\longrightarrow \Bbb C_{\varphi}$ is a point derivation at
$\varphi$.
 Let $a'',b''\in {\cal A}^{**}$
then there are nets $(a_\alpha)$ and $(b_\beta)$ in ${\cal A}$
such that converge respectively to $a''$ and $b''$ in the
$weak^*$- topology of ${\cal A}^{**}$. Then we have
\begin{eqnarray*}
({d_\varphi})''(a''b'')&=&weak^* lim_\alpha lim_\beta
{d_\varphi}(a_\alpha b_\beta) \\
&=&weak^* lim_\alpha lim_\beta {d_\varphi}(a_\alpha) \varphi
(b_\beta)+weak^* lim_\alpha lim_\beta \varphi (a_\alpha)
{d_\varphi}(b_\beta) \\
&=& ({d_\varphi})''(a'').\varphi ''(b'')+ \varphi''(a'').
({d_\varphi})''(b'').
\end{eqnarray*}
Thus $({d_\varphi})'':{ \cal A}^{**}\longrightarrow {\Bbb
C_{\varphi ''}}^{**}$ is a derivation. By Theorem 3.2 (ii), for
every  $a,b \in \cal A,$ we have
$$\langle d_\varphi (a),\varphi(b) \rangle+\langle d_\varphi
(b),\varphi(a) \rangle= \langle ({d_\varphi})''
(\hat{a}),{\varphi}''(\hat{b}) \rangle+ \langle ({d_\varphi})''
(\hat{b}),{\varphi}''(\hat{a}) \rangle=0.\hspace{2cm}
\hfill\blacksquare~$$

%

\end{document}